\documentclass[a4paper,12pt,leqno]{article}
\usepackage{latexsym}
\usepackage[all]{xy}

\usepackage{amssymb} 
\usepackage{amsmath} 

\def\versionnumber{ver.22}

\def\Z{{\mathbb{Z}}}
\def\R{{\mathbb{R}}}
\def\A{{\mathcal{A}}}

\DeclareMathOperator{\Der}{Der}

\DeclareMathOperator{\im}{im}

\numberwithin{equation}{section}

\newcommand{\bfm}{{\bf m}}
\newcommand{\bfa}{{\bf a}}
\newcommand{\bfb}{{\bf b}}

\newcommand{\bfh}{{\bf h}}
\newcommand{\bfz}{{\bf z}}
\newcommand{\bfc}{{\bf c}}

\newtheorem{theorem}{Theorem}[section]

\newtheorem{cor}[theorem]{Corollary}
\newtheorem{lemma}[theorem]{Lemma}
\newtheorem{define}[theorem]{Definition}

\title{
The freeness of
Shi-Catalan
arrangements
}
\author{
Takuro Abe
\footnote
{
Supported by JSPS Grants-in-Aid for Young Scientists
(B)
No. 21740014.
Department of Mechanical Engineering and Science,
Kyoto University,
Kyoto 606-8501, Japan.
email:abe.takuro.4c@kyoto-u.ac.jp
}
\,
and
Hiroaki Terao
\footnote
{
Supported by JSPS Grants-in-Aid, Scientific Research
(B) 
No. 21340001.
Department of Mathematics, Hokkaido University, 
Sapporo, Hokkaido 060-0810, Japan.
email:terao@math.sci.hokudai.ac.jp
}
(\versionnumber)}

\begin{document}

\maketitle

\begin{abstract}
Let 
$W$ be a finite Weyl group
and
$\A$ be the corresponding 
Weyl arrangement.
A deformation of $\A$
is an
affine
arrangement
which is obtained 
by adding 
to
each hyperplane
$H\in\A$
several 
parallel
translations of
 $H$ 
by the positive root
(and its integer multiples)
perpendicular to $H$.
We say that a deformation 
is $W$-equivariant if 
the number of parallel hyperplanes 
of each hyperplane 
$H\in \A$ depends only on 
the $W$-orbit of $H$.  
We prove that the conings of 
the $W$-equivariant deformations
are free arrangements
under a Shi-Catalan condition
and give a formula 
for the number of chambers. 
This
generalizes Yoshinaga's theorem 
conjectured by Edelman-Reiner.
%
\end{abstract}

\section{Introduction}
Let $V$ be an $\ell$-dimensional real vector space with 
an inner product $I: V \times V \rightarrow \R$.
Let
$S:=\mbox{Sym}(V^*)$ be the symmetric algebra of
the dual space $V^{*}$, 
$F$ the quotient field of $S$, and
$\Der_\R(S)$ the $S$-module of $\R$-linear derivations of $S$ to
itself. 
For a 
finite
{\bf
Weyl group} $W$, let us fix 
a positive system $\Phi_+$ with respect to $W$. 
For $\alpha\in\Phi_{+}$ define
\[
H_{\alpha} := \{v\in V \mid \alpha(v) = 0\} = \ker(\alpha),
\,\,\,
H_{\alpha, k} := \{v\in V \mid \alpha(v) = k\}
\,\,\,\,(k\in\Z).
\]
Then $
H_{\alpha, k}
$ 
is a parallel translation of $H_{\alpha} $. 
Then
$
\A:=\{ H_{\alpha}   \mid \alpha \in \Phi_+\}
$
is the {\bf
Weyl arrangement
}
 corresponding to $W$. 
A function $\bfm:\A \rightarrow \Z_{\ge 0}$ is called
a {\bf
multiplicity}.
For two multiplicities $\bfa, \bfb$,
 define
a {\bf deformation}
 $\A^{[-\bfa,\bfb]}$ of
 $\A$ by
$$
\A^ {[-\bfa,\bfb]}:=\{H_{\alpha, k} \mid -\bfa(H_{\alpha}) \le k \le 
\bfb(H_{\alpha})
,\ 
k\in \Z
,\
\alpha \in \Phi_+\},
$$
which
is an affine arrangement.
For basic concepts in the arrangement theory, consult
\cite{OT}.

\begin{define}
A multiplicity $\bfa$ is said to be {\bf $W$-equivariant}
if
$
\bfa(H)=\bfa(w H)
$ 
for   
every $H \in \A$ and
$w \in W$.
We say that $\A^{[-\bfa, \bfb]} $ is a {\bf Shi-Catalan} arrangement
if 
$\bfa$ and $\bfb$ are both 
{$W$-equivariant}
and 
$\im({\bf a}-{\bf b})\subseteq\{-1,0\}$.
\label{wid}
\end{define}

Suppose that $\A$ decomposes into $W$-orbits as
$
\A = \A_{1} \cup \dots \cup \A_{k}. 
$
%
Note that every $W$-orbit is not irreducible
because the type $B_{\ell}$-arrangement decomposes into
$A_{1}^{\ell}  $ and $D_{\ell} $.   
Identify a
$W$-equivariant multiplicity
$\bfm$ with
a
$k$-dimensional vector
$(m_{1}, \dots, m_{k})$ when 
$\bfm(H) = m_{j} \
(H\in \A_{j})$. 
Let $d_{j}^{(1)}, \dots, d_{j}^{(\ell)}$ 
be the exponents of $\A_{j}$
with
$d_{j}^{(1)}\le \dots\le d_{j}^{(\ell)}$
for
$1\le j\le k$. 
Let $h_{j} := d_{j}^{(\ell)} + 1$, which
is equal to the {\bf
Coxeter number}
 of $\A_{j} $ 
when $\A_{j} $ is irreducible.
Note that $\A_{j} $ is not irreducible only when
$\A_{j} $ is of the type $A_{1}^{\ell}$.
In this case $h_{j} = 2$
(see \cite{AT09, ATW}). 
Define
a
$k$-dimensional vector
$
\bfh := (h_{1} , \dots, h_{k}).
$ 
The following theorem is the main result of this article.

%
%

\begin{theorem}
\label{main} 
If $\A^{[-\bfa,\bfb]}$ is a Shi-Catalan
arrangement,
then
its coning
${\bf c} \left(\A^{[-\bfa,\bfb]}\right)$ is free.
\end{theorem} 

We have the following Corollary 
by Ziegler's theorem
(Theorem \ref{ziegler})
and Theorem \ref{ATWexp} 
by Wakamiko and the authors.

\begin{cor}
\label{mainexp} 
Suppose that $W$ is irreducible
and that $\A^{[-\bfa, \bfb]}$ is Shi-Catalan.
Then

(i) if $W$ is of the type $G_{2}$, $\bfa = \bfb=(b_{1}, b_{2})$,
and
$b_{1} + b_{2} $ is an odd integer,
then
the exponents of ${\bf c} \left(\A^{[-\bfb,\bfb]}\right)$
are
 given by
\[
1, 
2 + ({\bf b} \cdot {\bf h}),
4 + ({\bf b} \cdot {\bf h}).
\]

(ii) For all the other cases,
the exponents of ${\bf c}\left(\A^{[-\bfa,\bfb]}\right)$
are 
 given by
\[
1, 
m^{(1)} + ({\bf b} \cdot {\bf h}),
\dots
,
m^{(\ell)} + ({\bf b} \cdot {\bf h})
.
\]
Here the dot $\cdot$ stands for the ordinary inner product
of vectors
and
$m^{(1)}, \dots , m^{(\ell)}$ are
the exponents of
$
({\bf a}-{\bf b})^{-1} (0)$,
which is a union of $W$-orbits of $\A$.
\end{cor}

\noindent
{\em
Remark.} 
For any $W$ which may not be irreducible, 
the formula for the exponents is easily obtained from Corollary
\ref{mainexp}. 

\bigskip

%

We have the following formula for the number of chambers
by the factorization theorem
(Theorem \ref{factorization})
 in
\cite{Tf}
and Zaslavsky's theorem 
(Theorem \ref{zaslavsky})
in 
\cite{Zas}:

\begin{cor}
Suppose that $W$ is irreducible
and that $\A^{[-\bfa, \bfb]}$ is Shi-Catalan.

(i) If $W$ is of the type $G_{2}$, $\bfa = \bfb=(b_{1}, b_{2})$,
and
$b_{1} + b_{2} $ is an odd integer,
then
the number of chambers of $\A^{[-\bfb,\bfb]}$
is equal to
\[
(2 + 3b_{1} +3b_{2})(4 + 3b_{1} +3b_{2}).
\]

(ii) For all the other cases,
the number of chambers of $\A^{[-\bfa,\bfb]}$
is equal to
\[
\prod_{j=1}^{\ell}  
(m^{(j)} + ({\bf b} \cdot {\bf h}))
.
\]
Here 
$m^{(1)}, \dots , m^{(\ell)}$ are
the exponents of
$
({\bf a}-{\bf b})^{-1} (0)$.
\label{chamber}
\end{cor}

In particular, suppose that
the $W$-equivariant multiplicities
  $\bfa$ and $\bfb$ in Theorem \ref{main}
 are both constant. Then
Theorem \ref{main} (ii)
and Corollary
\ref{chamber} (ii)
are
Yoshinaga's results in
\cite{Y04} 
which were
conjectured by Edelman and Reiner
in \cite{ER}.
In this case,
the deformation $\A^{[-\bfa, \bfb]}$
is called an {\bf
extended Catalan arrangement}
when $\bfa=\bfb$ and
is called an 
{\bf
extended Shi arrangement}
when 
$\bfb=\bfa + 1$. 
Then
$
({\bf a}-{\bf b})^{-1} (0)$
is 
equal to 
$\A$ 
if  
$\bfa=\bfb$
and is the empty arrangement
if 
$\bfa\neq\bfb$.
The Shi-Catalan arrangement generalizes 
these two types of arrangements.

%

In Theorem \ref{mainexp}, the case
(i)
is the unique
obstruction for the formula in (ii)
to become a
blanket formula covering all the cases.
Recall that
both of $\A_{1} $ and $\A_{2} $ are of the
type $A_{2}$
in the $W$-orbit decomposition
when $\A$ is of the type 
$G_{2}$.
The reason of this unique
exception stems from
the fact that one cannot choose a 
$G_{2}$-invariant 
primitive derivation of $\A_{1}$ or $\A_{2}$
as shown in \cite{ATW}. 

\smallskip

The organization of this article is as follows. 
In Section 2, we review and summerize basic concepts 
and their properties.
In Section 3,
we reduce  
Theorem \ref{main}
to the irreducible
2-dimensional cases.
Our main tools are
Yoshinaga's freeness criterion
in \cite{Y04}
and
the freeness of $W$-equivariant 
Weyl multiarrangements in \cite{ATW}.
In Section 4, we complete the proof
of Theorem \ref{main} 
by verifying
the 2-dimensional cases
thanks to 
\cite{ER, ATH}
for $A_{2} $,
\cite{A09} for 
$B_{2} $, 
and
the addition-deletion theorem 
in
\cite{T}.


\section{Basic concepts and their properties}

An {\bf affine} arrangement of hyperplanes 
is a finite collection of
affine hyperplanes in $V$. 
If every hyperplane 
$H\in\A$ goes through the origin, then
$\A$ is called to be {\bf
central}.
When $\A$ is central, 
for each $H\in \A$ choose $\alpha_{H} \in V^{*} $ 
with $\ker (\alpha_{H} ) = H.$
Let $x_{1} , \dots, x_{\ell} $ be
a basis for the dual vector space $V^{*} $ of $V$. 
Let $\A$ be an affine arrangement in $V$
and $Q\in \R[x_{1} , \dots, x_{\ell}]$ 
be a defining polynomial for $\A$.  
Let $x_{0} $ be a new variable.
Let ${\bf c}Q\in
\R[x_{0}, x_{1} , \dots, x_{\ell}]$ 
be a homogeneous polynomial defined by
\[
{\bf c}Q := x_{0}^{1+\deg Q} 
Q(x_{1}/x_{0},
x_{2}/x_{0},
\dots,
x_{\ell}/x_{0}). 
\]
The {\bf
coning} ${\bf c}\A$ is a central arrangement
in ${\bf c}V :=
\R \oplus V$ defined by ${\bf c}Q$. 
Let $\overline{H_{\infty}} $ be the hyperplane 
in ${\bf c}V$ defined by 
$x_{0}=0 $.
For $H\in \A$ with
$H = \{\alpha = k\}$
\,\,$(\alpha\in V^{*}, k\in \R)$, 
let $\bfc H $ be the hyperplane in
${\bf c}V$ defined by
$
\alpha - k x_{0} = 0.
$  
If $\A = 
\{H_{1}, \dots, H_{n}  \}
$,
 then
${\bf c}\A
=
\{\overline{H_{\infty}}, \bfc{H_{1}}, \dots, \bfc{H_{n}}  \}.
$ 
For $Y\in L(\A)$ with 
$Y = H_{1} \cap \dots\cap H_{k}$, 
define $\bfc Y := \bfc H_{1} \cap \dots\cap
\bfc H_{k}
 \in L(\bfc \A)$. 

Let $\pi(\A, t)$ denote the {\bf Poincar\'e polynomial} 
\cite{OT}
of $\A$.
Since
$V$ is a real vector space,
each connected component of the complement
$V\setminus \cup_{H\in\A} H$ 
is called a {\bf chamber} of $\A$.
Recall 

\begin{theorem}
(Zaslavsky
\cite{Zas})
The number of chambers of $\A$ is equal to
$\pi(\A, 1)$. 
\label{zaslavsky} 
\end{theorem}

In the rest of this section, suppose that $\A$ is central.
Let $\Der_{S} $   
be the $S$-module of $\R$-linear derivations from $S$ to
itself.
Recall the derivation module
\[
D(\A)
=
\{
\theta\in\Der_{S}
\mid
\theta(\alpha_{H} ) \in \alpha_{H} S
{\rm ~for~all~}
H\in\A
\}
\]
over $S$.
We say that $\A$ is a {\bf
free arrangement} 
if $D(\A)$ is a
free
 $S$-module. 
%
When $\A$ is a free arrangement and $\theta_{1}, \dots, \theta_{\ell}  $
are a homogeneous $S$-basis for $D(\A)$, the integers
$\deg\theta_{1}, \dots, \deg\theta_{\ell}  $ are called the
{\bf exponents} of $\A$:
\[
\exp\A
=
(\deg\theta_{1}, \dots, \deg\theta_{\ell}).
\]
Every Weyl arrangement is a free arrangement
and their exponents are the same as the exponents of the
corresponding Weyl group by Saito (e.g., see 
\cite{S93}).
%

\begin{theorem}
(Factorization Theorem
\cite{T})
Suppose that $m^{(1)}, m^{(2)}, \dots , m^{(\ell)}$ 
are the exponents of a free arrangement $\A$.
Then $$\pi(\A, t) = \prod_{i=1}^{\ell} (1+ m^{(i)} t). $$ 
\label{factorization} 
\end{theorem}

For a multiplicity $\bfm : \A\rightarrow \Z_{\ge 0} $,
we call a pair $(\A, \bfm)$ a {\bf multiarrangement}.
Define 
$$
D(\A, \bfm)
:=
\{
\theta\in\Der_{S} ~|~ \theta(\alpha_{H} )\in \alpha_{H}^{\bfm(H)} S
\text{~for all~} H\in\A 
\}$$ 
which is a submodule of $\Der_{S} $.
We say that the multiarrangement $(\A, \bfm)$ is {\bf free}
if the $S$-module $D(\A, \bfm)$ is a free $S$-module.  
When $(\A, \bfm)$ is free, 
the exponenets of $(\A, \bfm)$, denoted by $\exp(\A, \bfm)$, 
are defined by
the degrees of a homogeneous $S$-basis for $D(\A, \bfm)$. 
Note that $\exp(\A, {\bf 1})=\exp\A$.

 For a given arrangement $\A$ and a fixed
hyperplane $H_{0}\in\A $,
define a multiarrangement
$(\A'', \bfz)$, which we call the {\bf Ziegler restriction}
\cite{Z}, by
\[
\A'' :=
\{
H_{0} \cap K ~|~ K\in \A':=\A\setminus\{H_{0} \}\},
\,\,
\bfz (X):= |\{
K\in\A'~|~ X=K\cap H_{0} 
\}|
\]
where $\A''$ is an arrangement living in $H_{0} $ 
and
$X\in\A''$. 
For any $Y\in L(\A)$ 
define the {\bf
localization} $\A_{Y} $ of $\A$ at $Y$ by
$
\A_{Y} :=\{H\in\A~|~Y\subseteq H\}.
$

\begin{theorem}
(Ziegler
\cite{Z})
If $\A$ is a free arrangement, then
\[
\exp(\A) = (1, d_{2}, \dots, d_{\ell})
\Leftrightarrow
\exp(\A'', \bfz) = (d_{2}, \dots, d_{\ell} ).
\]
\label{ziegler}
\end{theorem}

\begin{theorem}
\label{yoshinagastheorem}
(Yoshinaga's criterion
\cite{Y04})
Suppose $\ell>3$.
 For a central arrangement $\A$ and an 
arbitrary hyperplane
$H_{0}\in\A $, the following two conditions are
equivalent:

(1) $\A$ is a free arrangement,

(2) (2-i) the Ziegler restriction $(\A'', \bfz)$  
is free
and
(2-ii) 
$\A_{Y}$ is free for any $Y\in L(\A)\setminus \{{\bf 0}\}$
 such that $Y\subset H_{0} $. 
\end{theorem}

The following result is
by Table 4 in \cite{W} for (i)
and 
Theorem 1.3 in \cite{ATW} for (ii). 

\begin{theorem}
\label{ATWexp} 
%
Under the assumption of Corollary \ref{mainexp}, one has

(i) If $W$ is of the type $G_{2}$, $\bfa = \bfb=(b_{1}, b_{2})$,
and
$b_{1} + b_{2} $ is an odd integer,
$$\exp(\A(W), \bfa+\bfb + 1)=
(2 + 3b_{1} +3b_{2},
4 + 3b_{1} +3b_{2}).
$$

(ii) For all the other cases,
$$\exp(\A(W), \bfa + \bfb + 1)
=(
m^{(1)} + ({\bf b} \cdot {\bf h}),
\dots
,
m^{(\ell)} + ({\bf b} \cdot {\bf h})).
$$
In particular, the multiarrangements
$(\A(W), \bfa+\bfb+1)$ are free.
\end{theorem}


\section{The higher-dimensional cases}

Now we begin our proof of Theorem \ref{main} 
by an induction on the dimension $\ell$.
In this section we reduce the proof
to the two-dimensional cases by applying Yoshinaga's criterion
Theorem \ref{yoshinagastheorem}.
In
the next section we will complete the proof.

It suffices to verify the two conditions in Theorem 
\ref{yoshinagastheorem}.
The condition (2-i) follows from Theorem \ref{ATWexp}. 
Let us check the other condition (2-ii).
For this purpose we prove

\begin{lemma}

(i)
For any 
$X\in L({\bf c}
\left(
\A^{[-\bfa, \bfb]} 
\right) 
)
$ with 
$X\subseteq 
\overline{H_{\infty}}
$, there exists a unique 
$Y
\in
L(\A)
$ 
such that
$X
=
\bfc Y
\cap
\overline{H_{\infty}}
\in
L({\bf c}
\left(
\A^{[-\bfa, \bfb]} 
\right) 
)
$,  
and

\noindent
(ii)
$$
\left(
{\bf c}
\left(
\A^{[-\bfa, \bfb]} 
\right) 
\right)
_{X} 
=
{\bf c}
\left(
(\A_{Y})^{[-\bfa_{Y} , \bfb_{Y} ]} 
\right) 
=
{\bf c}
\left(
\A(W_{Y})^{[-\bfa_{Y} , \bfb_{Y} ]} 
\right), 
$$ 
where 
$\bfa_{Y} $
and
$\bfb_{Y} $
are the restrictions of
$\bfa$ and $\bfb$ respectively to 
$\A_{Y} $. 
\end{lemma} 

\noindent
{\em proof.}
(i)
 For $H\in \A^{[-\bfa, \bfb]} $, 
let $H^{(0)} \in \A $ be the parallel hyperplane 
of $H$ through the origin. 
Then $
\bfc H \cap \overline{H_{\infty}}
=
\bfc H^{(0)}  \cap \overline{H_{\infty}}.
$ 
There exist $H_{1}, \dots, H_{k}\in 
\A^{[-\bfa, \bfb]} 
$
such that
\[
X = 
\overline{H_{\infty}}
\cap \bfc H_{1} 
\cap \bfc H_{2} 
\cap
\dots
\cap \bfc H_{k}. 
\]
Define $Y = H^{(0)}_{1} 
\cap H^{(0)}_{2} 
\cap
\dots
\cap H^{(0)}_{k}
\in
L(\A).$ 

(ii)
The second equality follows from the fact that
the parabolic subgroup $W_{Y} $ is equal to
the group generated by the reflections with respect to
the Coxeter arrangement $\A_{Y} $.
We will prove the first equality.
Note
that
\begin{align*} 
\bfc H
\in
\left(
{\bf c}
\left(
\A^{[-\bfa, \bfb]} 
\right) 
\right)_{X} 
\Leftrightarrow
X=\bfc Y 
\cap 
\overline{H_{\infty}}
\subseteq
\bfc H
\Leftrightarrow
X=\bfc Y 
\cap 
\overline{H_{\infty}}
\subseteq
\bfc H^{(0)} \\
\Leftrightarrow
Y\subseteq H^{(0)}
\Leftrightarrow
H\in 
\left(
\A_{Y}
\right)^{[-\bfa_{Y}, \bfb_{Y}]}
\Leftrightarrow
\bfc H\in 
\bfc\left(\left(
\A_{Y}
\right)^{[-\bfa_{Y}, \bfb_{Y}]}
\right)
\end{align*} 
for $H\in \A^{[-\bfa, \bfb]}.$ 
This completes the proof.
$\square$ 

\medskip

Thanks to this lemma, the freeness of $\bfc\A^{[-\bfa, \bfb]} $ 
reduces into the freeness of 
${\bf c}
\left(
\A(W_{Y})^{[-\bfa_{Y} , \bfb_{Y} ]} 
\right)$
for each $Y\in L(\A(W))$ with $Y \neq \{\bf 0\}$. 
Note that $\A_{Y}$ is
a product of irreducible Weyl arrangements
of strictly lower ranks.
Also recall that
a product of free arrangements is also free.
Consequently the proof of Theorem \ref{main}
reduces to the irreducible
Weyl arrangements of rank two or lower. 
Note that the arrangement of the
type $A_{1} $ and the empty arrangement are obviously free.
Therefore, in the subequent section, we may assume that $\A$ is
an irreducible
two-dimensional Weyl arrangement in order to complete the proof
of
Theorem \ref{main}.

\section{The two-dimensional cases}

In this section we assume that $\A$ is either of the type
$A_{2} $, $B_{2} $ or
$G_{2} $.  

\medskip

{\bf
(A)
(the $A_{2} $ type)
}
As for $A_{2}$, Theorem \ref{main} is proved in \cite{ER} 
for 
the extended Catalan arrangements and in \cite{ATH}
for 
the extended Shi arrangements.

\medskip

{\bf
(B)
(the $B_{2} $ type)
}
As for $B_{2} $,
Theorem \ref{main} is proved in \cite{A09},
which we will review. 
Start from the free arrangement 
\begin{eqnarray*}
x&=&-sz,\ldots,sz\,\\
y&=&-sz,\ldots,sz\,\\
x\pm y&=&0,\\
z&=&0.
\end{eqnarray*}
where $z=0$ is the infinite hyperplane. The arrangement above is free 
with exponents $(1,2s+1,2s+3)$, the proof of which is an easy exercise and left to the reader. 
Define hyperplanes
\begin{align*}
H_{4a-3}: x-y=-az,
\,\,\,
H_{4a-2}:x+y=az,\\
H_{4a-1}:x-y=az,
\,\,\,
H_{4a}:x+y=-az\
\end{align*}
for $a=1,2,\ldots,t$.
Add $H_{1}, H_{2}, \dots, H_{4t}$ to the arrangement above in this order. 
Then the addition theorem completes the proof
in the case of the type $B_{2} $. 

  \medskip

{\bf
(G)
(the $G_{2} $ type)
}
Lastly we study the type $G_{2} $. 
Let $\A=\A(G_{2} )$ 
defined by 
$Q_1Q_2=0$ with 
\begin{eqnarray*}
Q_1&=&x(y-\displaystyle \frac{1}{\sqrt{3}}x)(y+\displaystyle \frac{1}{\sqrt{3}}x),\\
Q_2&=&y(y+\sqrt{3}x)(y-\sqrt{3}x).
\end{eqnarray*}
Put $\A_i:=\{Q_i=0\}$. Then 
$\A$ has an orbit decomposition
$$
\A=\A_1 \cup \A_2
$$
such that both $\A_1$ and $\A_2$ are of
the type $A_2$. 
Hence we have to verify the freenss of
the following 
four types of Shi-Catalan arrangements:

\medskip

{\bf
(G-i) (the Catalan-Catalan type)
}
the arrangement $\A[2s+1, 2t+1]$ 
is defined by
\begin{eqnarray*}
x&=&-s,\ldots,s\\
y-\displaystyle \frac{1}{\sqrt{3}}x&=&
-\displaystyle \frac{2}{\sqrt{3}}s,\ldots,\displaystyle \frac{2}{\sqrt{3}}s,\\
y+\displaystyle \frac{1}{\sqrt{3}}x&=&
-\displaystyle \frac{2}{\sqrt{3}}s,\ldots,\displaystyle \frac{2}{\sqrt{3}}s,\\
y&=&-\displaystyle \frac{1}{\sqrt{3}}t,\ldots,\displaystyle \frac{1}{\sqrt{3}}t,\\
y-\sqrt{3}x&=&-\displaystyle \frac{2}{\sqrt{3}}t,\ldots,\displaystyle \frac{2}{\sqrt{3}}t,\\
y+\sqrt{3}x&=&-\displaystyle \frac{2}{\sqrt{3}}t,\ldots,\displaystyle \frac{2}{\sqrt{3}}t
\end{eqnarray*}
with $s,t \in \Z_{\ge 0}$, 
which equals
$\A^{[-\bfa, \bfb]} $ where $\bfa=\bfb= (s, t)$. 

\medskip

{\bf
(G-ii) (the Shi-Catalan type)
}
The arrangement $\A[2s, 2t+1]$ 
is defined by
\begin{eqnarray*}
x&=&-s+1,\ldots,s\\
y-\displaystyle \frac{1}{\sqrt{3}}x&=&
-\displaystyle \frac{2}{\sqrt{3}}(s-1),\ldots,\displaystyle \frac{2}{\sqrt{3}}s,\\
y+\displaystyle \frac{1}{\sqrt{3}}x&=&
-\displaystyle \frac{2}{\sqrt{3}}(s-1),\ldots,\displaystyle \frac{2}{\sqrt{3}}s,\\
y&=&-\displaystyle \frac{1}{\sqrt{3}}t,\ldots,\displaystyle \frac{1}{\sqrt{3}}t,\\
y-\sqrt{3}x&=&-\displaystyle \frac{2}{\sqrt{3}}t,\ldots,\displaystyle \frac{2}{\sqrt{3}}t,\\
y+\sqrt{3}x&=&-\displaystyle \frac{2}{\sqrt{3}}t,\ldots,\displaystyle \frac{2}{\sqrt{3}}t
\end{eqnarray*}
with $s,t \in \Z_{\ge 0}$, 
 which equals
$\A^{[-\bfa, \bfb]} $ where $\bfa=(s-1, t)$ and $  \bfb= (s, t)$.

\medskip

{\bf
(G-iii)
(the Catalan-Shi type)
}
The arrangement
$\A[2s+1, 2t]$ 
is defined by
\begin{eqnarray*}
x&=&-s,\ldots,s\\
y-\displaystyle \frac{1}{\sqrt{3}}x&=&
-\displaystyle \frac{2}{\sqrt{3}}s,\ldots,\displaystyle \frac{2}{\sqrt{3}}s,\\
y+\displaystyle \frac{1}{\sqrt{3}}x&=&
-\displaystyle \frac{2}{\sqrt{3}}s,\ldots,\displaystyle \frac{2}{\sqrt{3}}s,\\
y&=&-\displaystyle \frac{1}{\sqrt{3}}(t-1),\ldots,\displaystyle \frac{1}{\sqrt{3}}t,\\
y-\sqrt{3}x&=&-\displaystyle \frac{2}{\sqrt{3}}(t-1),\ldots,\displaystyle \frac{2}{\sqrt{3}}t,\\
y+\sqrt{3}x&=&-\displaystyle \frac{2}{\sqrt{3}}(t-1),\ldots,\displaystyle \frac{2}{\sqrt{3}}t
\end{eqnarray*}
with $s,t \in \Z_{\ge 0}$, which equals
$\A^{[-\bfa, \bfb]} $ where $\bfa=(s, t-1)$ and $\bfb = (s, t)$.

\medskip

{\bf
(G-iv)
(the Shi-Shi type)
}
The arrangement $\A[2s, 2t]$ 
is defined by
\begin{eqnarray*}
x&=&-s+1,\ldots,s\\
y-\displaystyle \frac{1}{\sqrt{3}}x&=&
-\displaystyle \frac{2}{\sqrt{3}}(s-1),\ldots,\displaystyle \frac{2}{\sqrt{3}}s,\\
y+\displaystyle \frac{1}{\sqrt{3}}x&=&
-\displaystyle \frac{2}{\sqrt{3}}(s-1),\ldots,\displaystyle \frac{2}{\sqrt{3}}s,\\
y&=&-\displaystyle \frac{1}{\sqrt{3}}(t-1),\ldots,\displaystyle \frac{1}{\sqrt{3}}t,\\
y-\sqrt{3}x&=&-\displaystyle \frac{2}{\sqrt{3}}(t-1),\ldots,\displaystyle \frac{2}{\sqrt{3}}t,\\
y+\sqrt{3}x&=&-\displaystyle \frac{2}{\sqrt{3}}(t-1),\ldots,\displaystyle \frac{2}{\sqrt{3}}t
\end{eqnarray*}
with $s,t \in \Z_{\ge 0}$, 
which equals
$\A^{[-\bfa, \bfb]} $ where $\bfa=(s-1, t-1)$ and $\bfb = (s, t)$. 

%
%

We prove that these Shi-Catalan arrangements are all free by using the 
addition theorem. Since the counting of 
intersections
on a newly-added hyerplane
is easy, we just show 
the order of adding hyperplanes $H$ and the number of 
intersections $|H \cap \A|$ in the tables below.

\medskip

{\bf
(G-i) (G-iii) (when $s+t$ is odd)
}
First let us consider $\bfc\A[2s+1,2t+1]$ when  
$s+t$ is odd. We prove that $\bfc\A[2s+1,2t+1]$ is free with 
$$
\exp \bfc\A[2s+1,2t+1]=(1,A+2,A+4)
$$
where $A=3s+3t$. We use induction on $t \ge 0$. It is easy to check the freeness when 
$t=0$. Assume that $t \ge 1$. Then the addition table is as follows:
\begin{center}
\begin{tabular}{c|c|c}
\hline
added hyperplane & number of intersections & exponents\\
\hline
$y=\displaystyle \frac{t+1}{\sqrt{3}}z$ & $A+5$ & (1,A+3,A+4)\\
\hline
$y-\sqrt{3}x=\displaystyle \frac{2(t+1)}{\sqrt{3}}z$ & $A+5$ & (1,A+4,A+4)\\
\hline
$y+\sqrt{3}x=\displaystyle \frac{2(t+1)}{\sqrt{3}}z$ & $A+5$ & (1,A+4,A+5)\\
\hline
$y=-\displaystyle \frac{t+1}{\sqrt{3}}z$ & $A+5$ & (1,A+4,A+6)\\
\hline
$y-\sqrt{3}x=-\displaystyle \frac{2(t+1)}{\sqrt{3}}z$ & $A+5$ & (1,A+4,A+7)\\
\hline
$y+\sqrt{3}x=-\displaystyle \frac{2(t+1)}{\sqrt{3}}z$ & $A+5$ & (1,A+4,A+8)\\
\hline
$y=\displaystyle \frac{t+2}{\sqrt{3}}z$ & $A+9$ & (1,A+5,A+8)\\
\hline
$y-\sqrt{3}x=\displaystyle \frac{2(t+2)}{\sqrt{3}}z$ & $A+9$ & (1,A+6,A+8)\\
\hline
$y+\sqrt{3}x=\displaystyle \frac{2(t+2)}{\sqrt{3}}z$ & $A+9$ & (1,A+7,A+8)\\
\hline
$y=-\displaystyle \frac{t+2}{\sqrt{3}}z$ & $A+9$ & (1,A+8,A+8)\\
\hline
$y-\sqrt{3}x=-\displaystyle \frac{2(t+2)}{\sqrt{3}}z$ & $A+9$ & (1,A+8,A+9)\\
\hline
$y+\sqrt{3}x=-\displaystyle \frac{2(t+2)}{\sqrt{3}}z$ & $A+9$ & (1,A+8,A+10)\\
\hline
\end{tabular}
\end{center}

As a consequence, $\bfc\A[2s+1,2t+1]$ and $\bfc\A[2s+1,2t+2]$ are free with 
\begin{eqnarray*}
\exp \bfc\A[2s+1,2t+1]&=&(1,A+2,A+4),\\
\exp \bfc\A[2s+1,2t+2]&=&(1,A+4,A+5).
\end{eqnarray*}

\medskip

{\bf
(G-i) (G-iii) (when $s+t$ is even)
}
Next consider $G(2s+1,2t+1)$ when 
$s+t$ is even. Then the same table in the above shows that 
$\bfc\A[2s+1,2t+1]$ and $\bfc\A[2s+1,2t+2]$ are free with 
\begin{eqnarray*}
\exp \bfc\A[2s+1,2t+1]&=&(1,A+1,A+5),\\
\exp \bfc\A[2s+1,2t+2]&=&(1,A+4,A+5).
\end{eqnarray*}

\medskip

{\bf
(G-ii)(G-iv)
}
Next let us consider 
$\bfc\A[2s,2t+1]$ and $\bfc\A[2s,2t+2]$. We prove that they 
are free with 
\begin{eqnarray*}
\exp \bfc\A[2s,2t+1]&=&(1,A+1,A+2),\\
\exp \bfc\A[2s,2t+2]&=&(1,A+3,A+3).
\end{eqnarray*}
We begin with $\bfc\A[2s,2t+1]$. In this case 
the order of adding hyperplanes is important. The addition table 
is as follows:

\begin{center}
\begin{tabular}{c|c|c}
\hline
added hyperplane & number of intersections & exponents\\
\hline
$y=\displaystyle \frac{t+1}{\sqrt{3}}z$ & $A+3$ & $(1,A+2,A+2)$\\
\hline
$y+\sqrt{3}x=\displaystyle \frac{2}{\sqrt{3}}(t+1)z$ & $A+3 $& $(1,A+2,A+3)$\\
\hline
$y-\sqrt{3}x=\displaystyle \frac{2}{\sqrt{3}}(t+1)z$ &$ A+4$ & $(1,A+3,A+3)$\\
\hline
$y-\sqrt{3}x=-\displaystyle \frac{2}{\sqrt{3}}(t+1)z $& $A+4 $&$ (1,A+3,A+4)$\\
\hline
$y+\sqrt{3}x=-\displaystyle \frac{2}{\sqrt{3}}(t+1)z $& $A+5$ & $(1,A+4,A+4)$\\
\hline
$y=-\displaystyle \frac{t+1}{\sqrt{3}}z$ & $A+5 $& $(1,A+4,A+5)$	\\
\hline
\end{tabular}
\end{center}
Hence $\bfc\A[2s,2t+1]$ and $\bfc\A[2s,2t+2]$ 
are free with 
\begin{eqnarray*}
\exp \bfc\A[2s,2t+1]&=&(1,A+1,A+2),\\
\exp \bfc\A[2s,2t+2]&=&(1,A+3,A+3).
\end{eqnarray*}

The above tables show that each Shi-Catalan arrangements of the type $G_{2} $ 
is free.
%
%
Thus we complete the proof
of
Theorem \ref{main}. 


%
%
%
%

 \vspace{5mm}

\end{document}